\documentclass[11pt]{article}
\usepackage{amsmath,fullpage,amssymb,amsthm,hyperref}
\usepackage{tikz}

\def\M{\mathcal M}
\def\B{\mathcal B}
\def\P{\mathcal P}

\DeclareMathOperator\st{st}
\DeclareMathOperator\hs{hs}
\DeclareMathOperator\uhs{uhs}
\DeclareMathOperator\fd{fd}
\DeclareMathOperator\xfd{xfd}
\DeclareMathOperator\hfc{hfc}
\DeclareMathOperator\dr{dr}
\DeclareMathOperator\df{df}
\DeclareMathOperator\uc{ch_1}
\DeclareMathOperator\ch{ch}
\DeclareMathOperator\iuc{ich_1}
\DeclareMathOperator\uhu{uhu}
\DeclareMathOperator\hsp{hsp}
\DeclareMathOperator\hsv{hsv}
\DeclareMathOperator\lch{lch}
\DeclareMathOperator\lhs{lhs}
\DeclareMathOperator\stair{stair}
\DeclareMathOperator\oh{oh}
\DeclareMathOperator\eh{eh}
\DeclareMathOperator\area{area}
\def\WA{{\rm WA}}
\def\SA{{\rm SA}}

\def\U{-- ++(1,1) circle(1.2pt)}
\def\D{-- ++(1,-1) circle(1.2pt)}
\def\H{-- ++(1,0) circle(1.2pt)}
\def\u{-- ++(.2,.2) circle(0.2pt)}
\def\d{-- ++(.2,-.2) circle(0.2pt)}
\def\h{-- ++(.2,0) circle(0.2pt)}
\newcommand{\motzkin}[1]{\draw #1 circle(0.2pt) \u\u\h\d\h\d\h\u\h\u\h\h\d\h\u\u\h\d\d\d -- ++(-4,0);}
\def\N{-- ++(0,1) circle(1.2pt)}
\def\S{-- ++(0,-1) circle(1.2pt)}

\title{Statistics on bargraphs viewed as cornerless Motzkin paths}

\author{Emeric Deutsch\thanks{NYU Tandon School of Engineering, Brooklyn, NY 11201, USA.} \and Sergi Elizalde\thanks{Department of Mathematics, Dartmouth College, Hanover, NH 03755, USA. {\tt sergi.elizalde@dartmouth.edu}. Partially supported by grant \#280575 from the Simons Foundation and by grant H98230-14-1-0125 from the NSA.}
}

\begin{document}

\maketitle

\begin{abstract}
A bargraph is a self-avoiding lattice path with steps $U=(0,1)$, $H=(1,0)$ and $D=(0,-1)$ that starts at the origin and ends on the $x$-axis, and stays strictly above the $x$-axis everywhere except at the endpoints. 
Bargraphs have been studied as a special class of convex polyominoes, and enumerated using the so-called wasp-waist decomposition of Bousquet-M\'elou and Rechnitzer. In this paper we note that there is a trivial bijection between bargraphs and Motzkin paths without peaks or valleys. This allows us to use the recursive structure of Motzkin paths to enumerate bargraphs with respect to several statistics, finding simpler derivations of known results and obtaining many new ones.  We also count symmetric bargraphs and alternating bargraphs. In some cases we construct statistic-preserving bijections between different combinatorial objects, proving some identities that we encounter along the way.
\end{abstract}



\section{Introduction}

Bargraphs have appeared in the literature with different names, from skylines~\cite{Ger} to wall polyominoes~\cite{Fer}. Their enumeration was addressed by Bousquet-M\'elou and Rechnitzer~\cite{BMR}, Prellberg and Brak~\cite{PB}, and Fereti\'c~\cite{Fer}. These papers use the so-called {\em wasp-waist} decomposition or some variation of it to find a generating function with two variables $x$ and $y$ keeping track of the number of $H$ and $U$ steps, respectively.
The same decomposition has later been exploited by Blecher, Brennan and Knopfmacher to obtain refined enumerations with respect to statistics such as peaks~\cite{BBK_peaks}, levels~\cite{BBK_levels}, walls~\cite{BBK_walls}, descents and area~\cite{BBK_parameters}.
 Bargraphs are used to represent histograms and they also have connections to statistical physics, where they are used to model polymers.

There is a trivial bijection between bargraphs and Motzkin paths with no peaks or valleys, which surprisingly seems not to have been noticed before.
In Section~\ref{sec:bij} we present this bijection and we use it to describe our basic recursive decomposition. In Section~\ref{sec:stat} we find the generating functions for bargraphs with respect to several statistics, 
most of which are new, but in some cases we also present simpler derivations for statistics that have been considered before. In cases where we obtain enumeration sequences that coincide with each other or with sequences that have appeared in the literature, we construct bijections between the relevant combinatorial objects.
Finally, in Section~\ref{sec:subsets} we find the generating functions for symmetric bargraphs and for alternating bargraphs.

\section{Bijection to Motzkin paths}\label{sec:bij}

\begin{figure}[htb]
\centering
    \begin{tikzpicture}[scale=0.5]
     \draw(24,0) circle(1.2pt)  \N\N\N\H\S\H\S\H\N\H\N\H\H\S\H\N\N\H\S\S\S\S; 
      \draw[thin,dotted](24,0)--(32,0);
      \draw(22,1) node {$\rightarrow$};
      \draw(0,0) circle (1.2pt)  \U\U\H\D\H\D\H\U\H\U\H\H\D\H\U\U\H\D\D\D;  
       \draw[thin,dotted](0,0)--(20,0);
    \end{tikzpicture}
   \caption{A cornerless Motzkin path and the corresponding bargraph obtained by applying $\Delta$.}\label{fig:Delta}
\end{figure}
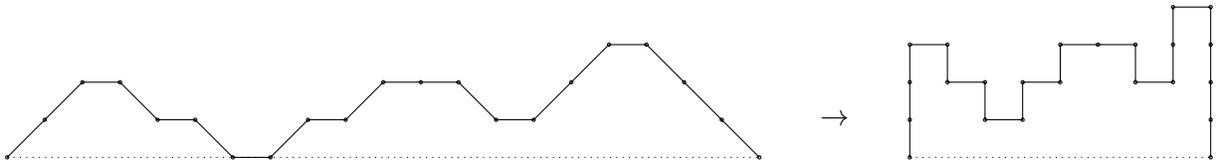

A bargraph is a lattice path with steps $U=(0,1)$, $H=(1,0)$ and $D=(0,-1)$ that starts at the origin and ends on the $x$-axis, stays strictly above the $x$-axis everywhere except at the endpoints, and has no pair of consecutive steps of the form $UD$ or $DU$. Sometimes it is convenient to identify a bargraph with the correspodning word on the alphabet $\{U,D,H\}$.
The {\em semiperimeter} of a bargraph is the number of $U$ steps plus the number of $H$ steps (this terminology makes sense when considering the closed polyomino obtained by connecting the two endpoints of the path with a horizontal segment). A bargraph of semiperimeter $15$ is drawn on the right of Figure~\ref{fig:Delta}. 
Let $\B$ denote the set of bargraphs, and let $\B_n$ denote those of semiperimeter $n$. Let $B=B(x,y)$ be the generating function for $\B$ where $x$ marks the number of $H$ steps and $y$ marks the number of $U$ steps. Note that the coefficient of $z^n$ in $B(z,z)$ is $|\B_n|$. We do not consider the trivial bargraphs with no steps, so the smallest possible semiperimeter is $2$.

Recall that a Motzkin path is a lattice path with up steps $(1,1)$, horizontal steps $(1,0)$ and down steps $(1,-1)$ that starts at the origin, ends on the $x$-axis, and never goes below the $x$-axis. A peak in a Motzkin path is an up step followed by a down step, and a valley is a down step followed by an up step. 

Let $\M$ denote the set of Motzkin paths without peaks or valleys, and let $\M_n$ denote the subset of those where the number of up steps plus the number of horizontal steps equals $n$. Elements of $\M$ will be called  {\em cornerless} Motzkin paths.
Let $M=M(x,y)$ be the generating function for $\M$ where $x$ marks the number of $H$ steps and $y$ marks the number of $U$ steps.
The coefficient of $z^n$ in $M(z,z)$ is $|\M_n|$. Let $\epsilon$ denote the empty Motzkin path, so that $\M_0=\{\epsilon\}$. To make our notation less clunky, we will sometimes write $\epsilon$ instead of $\{\epsilon\}$.

For $n\ge2$, there is an obvious bijection between $\M_{n-1}$ and $\B_n$ that appears not to have been noticed in the literature. 
Given a path in $\M$, insert an up step at the beginning and a down step at the end (resulting in an {\em elevated} Motzkin path), and then turn all the up steps $(1,1)$ into $U=(1,0)$, all the down steps $(1,-1)$ into $D=(0,-1)$, and leave the horizontal steps $(1,0)$ unchanged as $H=(1,0)$. Denote by $\Delta:\M_{n-1}\to\B_n$ the resulting bijection. An example is shown in Figure~\ref{fig:Delta}.

With some abuse of notation, we will denote up steps and down steps in a Motzkin path by $U$ and $D$ respectively, just like the steps in bargraphs that map to them via $\Delta$.

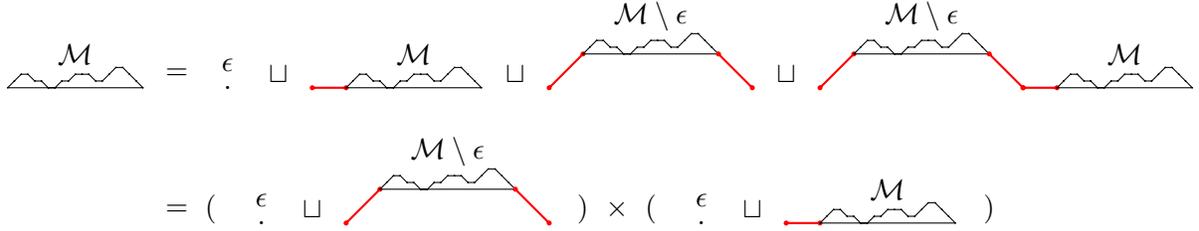
\begin{figure}[htb]
\centering
    \begin{tikzpicture}[scale=0.45]
     \motzkin{(-1,0)}
     \draw (1,0.4) node[above] {$\M$};
     \draw (4,0.4) node {$=$};
     \draw (5.5,0) node {$.$}; \draw (5.5,0.2) node[above] {$\epsilon$};
     \draw (7,0.4) node {$\sqcup$};
     \draw[thick,red] (8,0) circle(1.2pt) \H; \motzkin{(9,0)}; \draw (11,0.4) node[above] {$\M$};
     \draw (14,0.4) node {$\sqcup$};
      \draw[thick,red] (15,0) circle(1.2pt) \U; \motzkin{(16,1)}; \draw (18,1.4) node[above] {$\M\setminus \epsilon$}; \draw[thick,red] (20,1) circle(1.2pt) \D;
      \draw (22,0.4) node {$\sqcup$};
      \draw[thick,red] (23,0) circle(1.2pt) \U; \motzkin{(24,1)}; \draw (26,1.4) node[above] {$\M\setminus\epsilon$}; \draw[thick,red] (28,1) circle(1.2pt) \D;  \draw[thick,red] (29,0) circle(1.2pt) \H; \motzkin{(30,0)};  \draw (32,0.4) node[above] {$\M$};
      
       \begin{scope}[shift={(0,-4)}]
     \draw (4,0.4) node {$=$};
     \draw (5,0.4) node {$($};    
     \draw (6.5,0) node {$.$}; \draw (6.5,0.2) node[above] {$\epsilon$};
     \draw (8,0.4) node {$\sqcup$};
   \draw[thick,red] (9,0) circle(1.2pt) \U; \motzkin{(10,1)}; \draw (12,1.4) node[above] {$\M\setminus\epsilon$}; \draw[thick,red] (14,1) circle(1.2pt) \D;
     \draw (16,0.4) node {$)$};
     \draw (17,0.4) node {$\times$};
      \draw (18,0.4) node {$($};    
 \draw (19.5,0) node {$.$}; \draw (19.5,0.2) node[above] {$\epsilon$};
 \draw (21,0.4) node {$\sqcup$};
     \draw[thick,red] (22,0) circle(1.2pt) \H; \motzkin{(23,0)}; \draw (25,0.4) node[above] {$\M$};
          \draw (28,0.4) node {$)$};
  \end{scope}
      
    \end{tikzpicture}
   \caption{The decomposition of Motzkin paths without peaks or valleys.}\label{fig:decomposition}
\end{figure}

Cornerless Motzkin paths admit a very simple decomposition that parallels the standard recursive decomposition of Motzkin paths. Namely, every nonempty path in $\M$ can be decomposed uniquely as $HA$, $UA'D$, or $UA'DHA$, where $A$ is an arbitrary path in $\M$ and $A'$ is an arbitrary nonempty path. Figure~\ref{fig:decomposition} illustrates these cases. It follows that $M=M(x,y)$ satisfies the equation
\begin{align}\nonumber
M&=1+xM+y(M-1)+y(M-1)xM\\ &=(1+y(M-1))(1+xM).
\label{eq:M}
\end{align}
Using that $B=y(M-1)$, which is a consequence of the bijection $\Delta$, we obtain
\begin{equation}\label{eq:B} xB^2-(1-x-y-xy)B+xy=0.
\end{equation}
Solving for $B$, we get
$$B=\frac{1-x-y-xy-\sqrt{(1-x-y-xy)^2-4x^2y}}{2x},$$
which agrees with \cite[p.~92]{BMR}, \cite[Eq.~(3.6)]{BMB}, and~\cite[Eq.~(8)]{Fer} (with $y$ replaced by $y^2$).
The generating function for bargraphs by semiperimeter is
$$B(z,z)=\frac{1-2z-z^2-\sqrt{(1-2z-z^2)^2-4z^3}}{2z},$$
and the beginning of its series expansion is $z^2+2z^3+5z^4+13z^5+35z^6+97z^7+275z^8+\cdots$, 
given by sequence A082582 from the OEIS~\cite{OEIS}.

We point out that whereas the wasp-waist decomposition considers five cases by splitting bargraphs in different ways, the decomposition in Figure~\ref{fig:decomposition} is significantly simpler. It can in fact be considered as having one single case by interpreting the factorization $M=(1+y(M-1))(1+xM)$ as follows: every element of $\M$ can be uniquely split into two pieces: an elevated non-empty path (this piece could be empty), and a path that is either empty or starts with an H.

Via the bijection $\Delta$, the statistics on $\B$ that we will consider correspond naturally to statistics on $\M$.
When it creates no confusion, we will use the same notation for both a statistic on $\B$ and the corresponding  statistic on $\M$. 
Given such a statistic $\st$, we let $M_{\st}=M_{\st}(t,x,y)$
and $B_{\st}=B_{\st}(t,x,y)$  denote the generating functions for $\M$ and $\B$, respectively, where $t$ marks $\st$, $x$ marks the number of $H$, and $y$ marks the number of $U$. We will also use the analogous notation with the letter $M$ replaced by any symbol denoting a subset of $\M$.
If $\st(\epsilon)=0$ and $\st(M)=\st(\Delta(M))$ for all $M\in\M$, then $B_{\st}=y(M_{\st}-1)$ using the bijection $\Delta$. For some statistics these conditions do not hold, but we find similar relations between $B_{\st}$ and $M_{\st}$. Substituting $x=z$ and $y=z$ in $B_{\st}$ we obtain the generating function for $\B$ by semiperimeter and the statistic $\st$. 

In the study of some statistics, it will be helpful to separate cornerless paths according to the first and the last step.
The generating function of paths that start with an $H$ is $xM$, and the same applies for paths that end with an $H$.
Thus, the generating function for paths that start with a $U$ is $M-1-xM$, and similarly for paths that end with a $D$.
The generating function for paths that start and end with an $H$ is $x+x^2M$. By subtracting it from the generating function of paths that start with an $H$, it follows that paths that start with an $H$ and end with a $D$ are counted by $xM-x-x^2M$, and by symmetry so are paths that start with a $U$ and end with an $H$. 
By subtracting this from the generating function for paths that start with a $U$, we obtain that the generating function for paths that start with a $U$ and end with a $D$ is
$$M-1-xM-(xM-x-x^2M)=(1-x)^2M+x-1.$$
Table~\ref{tab:startend} summarizes these generating functions.

\begin{table}[h]
\centering
\begin{tabular}{c||c|c|c}
& end with an $H$ & end with a $D$ & total \\ \hline \hline
start with an $H$ & $x+x^2M $ & $(x-x^2)M -x$ & $xM $\\ \hline
start with a $U$ &  $(x-x^2)M -x$ & $(1-x)^2M +x-1$ & $M -1-xM $ \\ \hline
total & $xM $ &$M -1-xM $ & $M -1$
\end{tabular}
\caption{The generating functions for nonempty cornerless Motzkin paths with given initial and final steps.}
\label{tab:startend}
\end{table}

\section{Statistics on bargraphs}\label{sec:stat}

In this section we use the bijection $\Delta$ and the decomposition of cornerless Motzkin paths in Figure~\ref{fig:decomposition} to find the generating functions for bargraphs with respect to several statistics, while also keeping track of the number of $U$ and $H$ steps.
For each statistic $\st$ we obtain an explicit formula for $B_{\st}$. In some cases, for the sake of simplicity, we will only state the equation that it satisfies, which is quadratic for all the statistics in this paper.

The statistics that we consider are the height of the first column (Section \ref{hfc}); the number of double rises and double falls (\ref{dr});
the number of valleys or peaks of width $\ell$ along with the number of horizontal steps in such valleys and peaks (\ref{vl});
the number of corners of different types (\ref{cor}); the number of horizontal segments (\ref{hs}) and horizontal segments of length 1 (\ref{uhs}); the length of the first descent, together with a bijective proof of a simple recurrence satisfied by this statistic, and the $x$-coordinate of the first descent (\ref{fd}); the number of columns of height $h$ for any fixed $h$, and
the number of initial columns of height 1 (\ref{ch}); the least column height~(\ref{lch}) and the width of the leftmost horizontal segment (\ref{lhs}), together with a bijective proof of their equidistribution on bargraphs with fixed semiperimeter, and a bijection relating the latter statistic with the number of initial columns of height 1 (\ref{lhs}); the number of occurrences of $UHU$ (\ref{uhu}); the length of the initial staircase $UHUH\dots$ (\ref{stair}); the number of odd-height and even-height columns (\ref{oheh}); and the area~(\ref{area}).
 Table~\ref{tab:OEIS} lists the OEIS reference~\cite{OEIS} of the sequences that occur in the paper, most of which have been added as a result of this work.

\begin{table}[htb]
\centering
\begin{tabular}{l|l|l}
Sect. & Statistic or subset of bargraphs & OEIS sequences \\
\hline\hline
\ref{hfc} & height of the first column & A273342, A273343 \\ \hline
\ref{dr} & number of double rises & A273713, A273714 \\ \hline
 \ref{dr} &   number of double rises and double falls & A276066 \\ \hline
\ref{vl} & number of valleys of width 1 & A273721, A273722 \\ \hline
 \ref{vl} &   number of peaks of width 1  &  A273715, A273716 \\ \hline
 \ref{cor} & number of $\llcorner$ corners  &  A273717, A273718 \\ \hline
\ref{hs} & number of horizontal segments  &  A274486 \\ \hline
\ref{uhs} & number of horizontal segments of length $1$  &  A274491  A274492 \\ \hline
\ref{fd} & length of the first descent  &   A276067,   A276068 \\ \hline
\ref{fd} &  $x$-coordinate of the first descent  & A273897, A273898 \\ \hline
\ref{ch} & number of columns of height $1$  &  A273899, A273900 \\ \hline
\ref{ch} & number of initial columns of height $1$ &  A274490 \\ \hline
\ref{lch} & least column height & A274488 \\ \hline
\ref{lhs} & width of the leftmost horizontal segment & A274488 \\ \hline
\ref{uhu} & number of occurrences of $UHU$ &  A273896 \\ \hline
\ref{stair} & length of the initial staircase $UHUH\dots$ & A274494, A274495 \\ \hline
\ref{oheh} & number of odd-height and even-height columns & A273901,
A273902, A273903, A273904 \\ \hline
\ref{area} & area & A273346, A273347, A273348 \\ \hline
\ref{sec:sym} & symmetric bargraphs  & A273905 \\ \hline
\ref{sec:alt} & weakly alternating bargraphs  & A275448 \\ \hline
 \ref{sec:alt} & strictly alternating bargraphs  & A023342 \\
\end{tabular}
\caption{The OIES sequences~\cite{OEIS} corresponding to the statistics and subsets of bargraphs studied in this paper.}
\label{tab:OEIS}
\end{table}

We use the term {\em ascent} (resp. {\em descent}) to refer to a maximal consecutive sequence of $U$ (resp. $D$) steps.

\subsection{Height of the first column}\label{hfc}

Let $\hfc$ be the statistic ``height of the first column'' on $\B$. 
On $\M$, let $\hfc$ denote the statistic that measures the length of the initial ascent of the path, that is, the number of steps $U$ right at the beginning of the path.
 Then, using the second line of the decomposition in Figure~\ref{fig:decomposition}, we get
$$M_{\hfc}=(1+ty(M_{\hfc}-1))(1+xM),$$
where we can easily use the expression for $M$ coming from~\eqref{eq:M} and solve for $M_{\hfc}$.

If $A\in\M$ and $G=\Delta(M)$, then $\hfc(G)=\hfc(A)+1$. It follows that $B_{\hfc}=ty(M_{\hfc}-1)$, and
$B_{\hfc}$ satisfies the equation
$$(1 - t(1 -x + y + xy) + t^2 y) B_{\hfc}^2 - t(1 - y)(1 - x - ty - txy) B_{\hfc} + t^2 xy(1 - y) = 0.$$
In particular, the generating function for $\B$ with respect to semiperimeter and height of the first column is
$$B_{\hfc}(t,z,z)=
t \,\frac {1-2z-tz+{z}^{2}+t{z}^{3} -\left( 1-tz \right) \sqrt{(1-z)(1-3z-z^2-z^3)} }{2(1-t+t^2 z-t{z}^{2})}.
$$

\subsection{Number of double rises and number of double falls}\label{dr}

Let $\dr$ (resp. $\df$) be the statistic that counts pairs of adjacent $U$ steps (resp. $D$ steps) ---also called {\em double rises} (resp. {\em double falls})--- in bargraphs and in cornerless Motzkin paths. By symmetry, $\dr$ and $\df$ are equidistributed.
Using the second line of the decomposition in Figure~\ref{fig:decomposition}, we have
\begin{equation}\label{eq:Mdr}
M_{\dr}=\left[1+y(t(M_{\dr}-1-xM_{\dr})+xM_{\dr})\right](1+xM_{\dr}).
\end{equation}
This is because, in this decomposition, the number of double rises in a path is obtained by summing the number of double rises in each of the two pieces, except that if the first piece is the elevation of a path starting with a $U$, it creates an additional double rise.
The generating function with respect to $\dr$ for paths starting with an $H$ is $xM_{\dr}$, and for paths starting with a $U$ is $M_{\dr}-1-xM_{\dr}$.
Equation~\eqref{eq:Mdr} can be solved to find an expression for $M_{\dr}$.

For every $A\in\M$, $\dr(\Delta(A))$ equals $\dr(A)+1$ if $A$ starts with an $U$ step and it equals $\dr(A)$ otherwise. It follows that 
$B_{\dr}=y\left(t(M_{\dr}-1-xM_{\dr})+xM_{\dr}\right)$, and so it satisfies
$$x B_{\dr}^2  - (1 - x  - ty- xy) B_{\dr} + xy = 0.$$
In particular, the generating function for bargraphs with no double rises according to semiperimeter is
$$B_{\dr}(0,z,z)=\frac {1-z-z^2-\sqrt {1-2z-z^2-2z^3+z^4}}{2z},$$
which agrees with the generating function for RNA secondary structure numbers~\cite[Exer. 6.43]{EC2} with the exponents shifted by one. Recall that a {\em secondary structure} is a simple graph with vertices $\{1,2,\dots,n\}$ such that every $i$ is adjacent to at most one vertex $j$ with $|j-i|>1$, $\{i,i+1\}$ is an edge for all $i$, and there are no two edges $\{i,k\},\{j,l\}$ with $i<j<k<l$.

Next we give a bijection between bargraphs of semiperimeter $n+1$ with no double rises and secondary structures on $n$ vertices. Given such a bargraph, delete the initial $U$ and the final $D$. Since the bargraph has no double rises, the remaining steps can be written uniquely as a sequence of blocks $HU$, $H$ and $D$. Note that there are $n$ blocks, and that no $HU$ is immediately followed by a $D$. Each block will correspond to a vertex of the secondary structure, labeled increasingly from left to right.  For each $HU$ block, consider the first $D$ block to its right that is at the same height in the bargraph, and draw an edge between the two corresponding vertices in the secondary structure. Finally, draw edges between every pair of consecutive vertices. It is easy to see that this construction produces a secondary structure, and that it is a bijection.

We end this section by refining the argument behind Equation~\eqref{eq:Mdr} in order to give the joint distribution of $\dr$ and $\df$ on bargraphs and cornerless Motzkin paths. Let $B_{\dr,\df}$ and $M_{\dr,\df}$ be the corresponding generating functions where $t$ marks $\dr$ and $s$ marks $\df$. In the decomposition  in Figure~\ref{fig:decomposition}, the elevated piece of the path creates an additional double rise if it starts with a $U$, and an additional double fall if it ends with a $D$.
The generating function by $\dr$ and $\df$ of paths that start with an $H$ is $xM_{\dr,\df}$, since the initial $H$ step does not affect the statistics $\dr$ and $\df$, and the same holds for paths that end with an $H$. The generating function by $\dr$ and $\df$ of paths that start and end with an $H$ is $x+x^2M_{\dr,\df}$. Thus, by the same argument used to obtain Table~\ref{tab:startend}, but keeping track of $\dr$ and $\df$, it follows that paths that start with an $H$ and end with a $D$ are counted by $xM_{\dr,\df}-x-x^2M_{\dr,\df}$, and so are paths that start with a $U$ and end with an $H$. Similarly, the generating function for paths that start with a $U$ and end with a $D$ is
$(1-x)^2M_{\dr,\df}+x-1$.

Adding the contributions of the elevated path in the different cases, we get the equation
$$
M_{\dr,\df}=\left[1+y\left[x+x^2M_{\dr,\df}+(t+s)((x-x^2)M_{\dr,\df}-x)+ts((1-x)^2M_{\dr,\df}+x-1)\right]\right](1+xM_{\dr,\df}).
$$
Using that $B_{\dr,\df}=y\left[x+x^2M_{\dr,\df}+(t+s)((x-x^2)M_{\dr,\df}-x)+ts((1-x)^2M_{\dr,\df}+x-1)\right]$, we obtain the following equation for $B_{\dr,\df}$:
$$xB_{\dr,\df}^2-(1-x-(t+s)xy-tsy+tsxy)B_{\dr,\df}+xy=0.$$

\subsection{Number of valleys and peaks of width $\ell$}\label{vl}

For $\ell\ge1$, define a {\em valley} (resp. {\em peak}) of width $\ell$ in a bargraph or a cornerless Motzkin path to be a sequence of steps $DH^\ell U$ (resp. $UH^\ell D$). Let $v_\ell$ (resp. $p_\ell$) be statistic counting the number of valleys (resp. peaks) of width $\ell$.

To count valleys of width $\ell$, we use the first line in the decomposition in Figure~\ref{fig:decomposition} to obtain
$$M_{v_\ell}=1+x M_{v_\ell}+y(M_{v_\ell}-1)+xy(M_{v_\ell}-1)\left[M_{v_\ell}+(t-1)x^{\ell-1}(M_{v_\ell}-1-xM_{v_\ell})\right].$$
Indeed, an extra valley of width $\ell$ is created in the last case when the second part of the path starts with $H^{\ell-1}U$; these paths are counted by $x^{\ell-1}(M_{v_\ell}-1-xM_{v_\ell})$. Using that $B_{v_\ell}=y(M_{v_\ell}-1)$, we obtain the following equation for $B_{v_\ell}$:
$$(1 - (1 - t)(1 - x)x^{\ell-1})B_{v_\ell}^2 - (1-x-y-xy- (1 - t)x^{\ell+1}y) B_{v_\ell} + xy = 0.$$
Solving it and setting $x=z$ and $y=z$, we get
$$B_{v_\ell}(t,z,z)= \frac{1-2z-{z}^{2}+(t-1){z}^{\ell+2}-\sqrt{1-4z+2{z}^{2}+{z}^{4}
+2 (1-t) {z}^{\ell+2}(1+z^2)
+(1-t)^2 {z}^{2\ell+4}}}
{2(z-(1-t)(1-z)z^\ell)}.$$

Similarly, for peaks of width $\ell$,  using the second line in the decomposition in Figure~\ref{fig:decomposition}, we have
$$M_{p_\ell}=[1+y(M_{p_\ell}-1+(t-1)x^\ell)](1+x M_{p_\ell}),$$
since an extra peak of width $\ell$ is created when the path that is elevated is precisely $H^\ell$.
Using that $B_{p_\ell}=y(M_{p_\ell}-1+(t-1)x^\ell)$, we obtain the equation
$$x B_{p_\ell}^2 - (1-x-y-xy- (1 - t)x^{\ell+1}y ) B_{p_\ell} +  y(x - (1 - t)(1 - x)x^\ell)= 0.$$
It follows that
$$B_{p_\ell}(t,z,z)= \frac{1-2z-{z}^{2}+(t-1){z}^{\ell+2}-\sqrt{1-4z+2{z}^{2}+{z}^{4}
+2 (1-t) {z}^{\ell+2}(1+z^2)
+(1-t)^2 {z}^{2\ell+4}}}
{2z}.$$

Generating functions for several statistics related to peaks in bargraphs were found in~\cite{BBK_peaks} using the wasp-waist decomposition. The statistics include the total number of peaks, the number of horizontal steps in peaks, and the height of the first peak. All of them can be derived using the decomposition in Figure~\ref{fig:decomposition}. For example, to find the generating function $B_{p,\hsp}(t,s,x,y)$ for bargraphs where $t$ marks the number of peaks and $s$ marks the number of horizontal steps in peaks, we first note that the corresponding generating function $M_{p,\hsp}(t,s,x,y)$ for cornerless Motzkin paths satisfies
$$M_{p,\hsp}=\left(1+y\left(M_{p,\hsp}-1+\frac{tsx}{1-sx}-\frac{x}{1-x}\right)\right)(1+xM_{p,\hsp}).$$
This is because in the decomposition in Figure~\ref{fig:decomposition}, when the elevated path in the first piece is of the form $H^\ell$ for $\ell\ge1$, it creates a peak of width $\ell$, contributing $ts^\ell x^\ell$ to the generating function instead of $x^\ell$. Summing $ts^\ell x^\ell-x^\ell$ for $\ell\ge1$ gives the term $\frac{tsx}{1-sx}-\frac{x}{1-x}$. Solving for $M_{p,\hsp}$ and using that $$B_{p,\hsp}=y\left(M_{p,\hsp}-1+\frac{tsx}{1-sx}-\frac{x}{1-x}\right)$$ we obtain an expression for $B_{p,\hsp}$.

We can similarly find the generating function $B_{v,\hsv}(t,s,x,y)$ where $t$ marks the number of valleys and $s$ marks the number of horizontal steps in valleys. In this case, the corresponding generating function $M_{v,\hsv}(t,s,x,y)$ for
$\M$ satisfies
$$M_{v,\hsv}=1+xM_{v,\hsv}+y(M_{v,\hsv}-1)+xy(M_{v,\hsv}-1)\left[H+\left(\frac{ts}{1-sx}-\frac{1}{1-x}\right)(M_{v,\hsv}-1-xM_{v,\hsv})\right].$$
Indeed, in the fourth case of the decomposition in the first line of Figure~\ref{fig:decomposition}, an extra valley of width $\ell$ appears if the path to the right of the red $H$ step starts with $H^{\ell-1}U$ for $\ell\ge1$. Such paths thus contribute $ts^\ell x^{\ell-1}(M_{v,\hsv}-1-xM_{v,\hsv})$ instead of $x^{\ell-1}(M_{v,\hsv}-1-xM_{v,\hsv})$.
Solving for $M_{v,\hsv}$ and using that $B_{v,\hsv}=y(M_{v,\hsv}-1)$ we obtain an expression for $B_{v,\hsv}$.

\subsection{Number of corners}\label{cor}
Denote by $\llcorner$, $\lrcorner$, $\ulcorner$ and $\urcorner$ the statistics counting the number of occurrences of $DH$, $HU$, $UH$ and $HD$ in a bargraph, respectively. For every $G\in\B$, we have $\ulcorner(G)=1+\lrcorner(G)$ and $\urcorner(G)=1+\llcorner(G)$. Also, by symmetry, the statistics $\ulcorner$ and $\urcorner$ are equidistributed on bargraphs, and so are $\llcorner$ and $\lrcorner$. The statistic $\urcorner$ coincides with the number of descents, and similarly the statistic $\ulcorner$ equals the number of ascents.

We can find the distribution of these statistics using the decomposition for $\M$ in the first line of Figure~\ref{fig:decomposition}. Let $\llcorner$ and $\ulcorner$ denote the number of $DH$ and $UH$ in paths in $\M$, respectively, and let $M_{\llcorner,\ulcorner}=M_{\llcorner,\ulcorner}(t,s,x,y)$ be the generating function where $t$ and $s$ mark these two statistics, respectively. Then
\begin{equation}\label{eq:MDHUH}
M_{\llcorner,\ulcorner}=1+xM_{\llcorner,\ulcorner}+y(M_{\llcorner,\ulcorner}-1+(s-1)xM_{\llcorner,\ulcorner})(1+txM_{\llcorner,\ulcorner}),
\end{equation}
since an extra occurrence of $UH$ is created when the elevated piece of the decomposition starts with an $H$, and an extra occurrence of $DH$ is created in the fourth case of the decomposition.
The generating function for bargraphs where $s$ and $t$ mark the statistics $\llcorner$ and $\ulcorner$, respectively, is then $B_{\llcorner,\ulcorner}=y(M_{\llcorner,\ulcorner}-1+(s-1)xM_{\llcorner,\ulcorner})$. 
Eliminating $M_{\llcorner,\ulcorner}$ from this equation and~\eqref{eq:MDHUH}, we see that $B_{\llcorner,\ulcorner}$ satisfies
$$txB_{\llcorner,\ulcorner}^2 - (1 - x - y + xy - txy - sxy)B_{\llcorner,\ulcorner} + sxy = 0.$$

Setting $t=0$ and $s=1$ in the above equation, we see that the generating function for bargraphs with no occurrences of $DH$ is $\frac{xy}{1-x-y}$. These are {\em nondecreasing} bargraphs, namely those where the heights of the columns increase weakly from left to right. Their generating function can be easily explained combinatorially, since such a bargraph is determined by the left boundary, which consists of any sequence of $U$ and $H$ steps that starts with a $U$ and ends with an $H$. A related family is that of {\em increasing} bargraphs, which are those avoiding $DH$ and $HH$. Their generating function is $\frac{xy}{1-y-xy}$, since the left boundary has the additional restriction of not containing two consecutive $H$s.

Since the total number of corners in a bargraph is counted by $\llcorner+\lrcorner+\ulcorner+\urcorner=2(\llcorner+\ulcorner)$, the generating function 
$$
B_{\llcorner,\ulcorner}(u^2,u^2,x,y)=\frac {1-x-y+(1-2u^2)xy-
\sqrt { (1-y ) ( 1-2x-y+(2-4{u}^{2})xy+{x}^{2}-
(1-2u^2)^2{x}^{2}y ) }}{2{u}^{2}x}
$$
enumerates bargraphs with respect to the total number of corners, marked by $u$. 
Viewing the bargraph as a polyomino, the two additional bottom corners can be easily included by multiplying the generating function by $u^2$.

In~\cite{BBK_walls}, the authors find the distribution of the number of ascents (called {\em walls} in that paper) of any fixed length $r$.

\subsection{Number of horizontal segments} \label{hs}

Define a {\em horizontal segment} to be a maximal sequence of consecutive $H$ steps, and 
let $\hs$ be the statistic ``number of horizontal segments,'' both on $\B$ and on $\M$. 
Horizontal segments of length at least 2 are considered in~\cite{BBK_levels}, where they are called {\em levels}.

In any bargraph, the number of horizontal segments equals half of the total number of corners.
It follows that $\B_{\hs}=B_{\llcorner,\ulcorner}(t,t,x,y)$, and so
\begin{equation}
\label{eq:Bhs}
B_{\hs}=\frac {1-x-y+xy-2txy-\sqrt{(1-y)(1-2x-y+{x}^{2}+(2-4t)xy-(1-2t)^2{x}^{2}y)}}{2tx}.
\end{equation}

Alternatively, a direct derivation can be given as follows. 
Let $\M^H$ be the subset of $\M$ consisting of paths that start with an $H$. For $A\in\M$, he have that 
$\hs(HA)$ equals $\hs(A)$ if $A\in\M^H$ and it equals $\hs(A)+1$ otherwise.
Via the usual decomposition, we have
\begin{align*}
M^H_{\hs}&=x[t(M_{\hs}-M^H_{\hs})+M^H_{\hs}],\\
M_{\hs}&=1+M^H_{\hs}+y(M_{\hs}-1)(1+M^H_{\hs}).
\end{align*}
Eliminating $M^H_{\hs}$ from these equations and using that $B_{\hs}=y(M_{\hs}-1)$, we recover~\eqref{eq:Bhs}.

\subsection{Number of horizontal segments of length 1} \label{uhs}

Let $\uhs$ denote the statistic ``number of horizontal segments of length 1," also called unit horizontal segments, both on $\B$ and on $\M$. 
Let $\M^{(1)}$ be the subset of $\M$ consisting of paths that start with a horizonal segment of length 1 (note that these are the paths that start with $HU$, plus the length-one path $H$). Let $\M^{(2)}$ be the subset of $\M$ consisting of paths that start with $HH$. Note that the paths in $\M$ but not in $\M^{(1)}$ or $\M^{(2)}$ are precisely those that are empty or start with a $U$. 

For $A\in\M$, 
$$\uhs(HA)=\begin{cases} 
\uhs(A)-1 & \text{if $A\in\M^{(1)}$,} \\
\uhs(A) & \text{if $A\in\M^{(2)}$,} \\
\uhs(A)+1 & \text{otherwise.}
\end{cases}$$
Using the standard decomposition, it follows that
\begin{align*}
M^{(1)}_{\uhs}&=tx(M_{\uhs}-M^{(1)}_{\uhs}-M^{(2)}_{\uhs}),\\
M^{(2)}_{\uhs}&=x\left(\frac{M^{(1)}_{\uhs}}{t}+M^{(2)}_{\uhs}\right),\\
M_{\uhs}&=1+M^{(1)}_{\uhs}+M^{(2)}_{\uhs}+y(M_{\uhs}-1)(1+M^{(1)}_{\uhs}+M^{(2)}_{\uhs}).
\end{align*}
Solving this system of equations, we get an expression for $M_{\uhs}$. Using that $B_{\uhs}=y(M_{\uhs}-1)$, we obtain the equation
$$x(t + x - tx)B_{\uhs}^2 -(1-x-y+xy- 2xy(t + x - tx))B_{\uhs} + xy(t + x -tx) = 0.$$

\subsection{Length and abcissa of the first descent} \label{fd}
Let $\fd$ denote the length of the first descent in a bargraph. The corresponding statistic on $\M$, via $\Delta$, is the length of the first descent of the cornerless Motzkin path, with the caveat that we have to add one if this first descent occurs at the end of the path. We also denote this statistic by $\fd$.

Let  $\M^I$ be the subset of $\M$ consisting of Motzkin paths that have the first descent (equivalently, all their $D$ steps) at the end. These correspond to nondecreasing bargraphs, mentioned in Section~\ref{cor}. Clearly,
$$M^I_{\fd}=\frac{tx}{1-x-ty},$$
since such a path is determined by a sequence of $U$ and $H$ steps that ends with an $H$, and only the $U$ steps contribute to $\fd$, with the $t$ in the numerator accounting for the fact that we have to add one when the first descent occurs at the end.

For a path of the form $UADB\in\M$ (where $A\in\M$ is nonempty and $B\in\M$ is either empty or starts with an $H$), we have $\fd(UADB)=\fd(A)$ except when $A\in\M^I$ and $B$ is empty, in which case $\fd(UADB)=\fd(A)+1$. It follows that
$$M_{\fd}=t+xM_{\fd}+y(M_{\fd}-t+(t-1)M^I_{\fd})+y(M_{\fd}-t)xM.$$
We can now use the known expressions for $M^I_{\fd}$ and $M$ and solve for $M_{\fd}$. Since $B_{\fd}=y(M_{\fd}-t)$, it follows that
\begin{equation}\label{eq:Bfd}B_{\fd}=t(1-x-y)\frac{1 - x - y - xy-\sqrt{(1-y)((1-x)^2 - y(1 + x)^2)}}{2 x (1-x-ty)}.
\end{equation}

If we denote by $a_{n,k}$ the number of bargraphs of semiperimeter $n$ whose first descent has length $k$, so that 
$$B_{\fd}(t,z,z)=\sum_{k,n}a_{n,k}t^kz^n,$$
then Equation~\eqref{eq:Bfd} implies that
\begin{equation}\label{eq:Bfd2}(1-x-tz)\sum_{k,n}a_{n,k}t^kz^n=t g(z)
\end{equation}
for some function $g(z)$ that does not depend on $t$.
For any fixed $k,n\ge2$, extracting the coefficient of $t^kz^n$ on both sides of~\eqref{eq:Bfd2} gives
$$a_{n,k}-a_{n-1,k}-a_{n-1,k-1}=0.$$

This simple recurrence can also be proved directly by the following combinatorial argument. Let $A_{n,k}$ denote the set of bargraphs with semiperimeter $n$ and first descent of length $k$, so that $a_{n,k}=|A_{n,k}|$.
Fix $n,k\ge2$. Given a bargraph in $A_{n,k}$, the step before the first descent must be an~$H$. The step before that, which we denote by $s$, must be either an $H$ or a $U$. This splits the set $A_{n,k}$ into two disjoint sets $A^H_{n,k}$ and $A^U_{n,k}$, where the superindex indicates the step $s$. 
A bijection between $A^H_{n,k}$ and $A_{n-1,k}$ is obtained by removing the step $s=H$. 
A bijection between $A^U_{n,k}$ and $A_{n-1,k-1}$ is obtained by removing the step $s=U$ and a $D$ from the first descent of the path. For example, for $n=6$ and $k=3$,
the bargraph $UUUUHHHDDDD$ is mapped to $UUUHHDDD$, and $UUUUHDDDHD$ is mapped to $UUUHDDHD$.
It follows that 
$$a_{n,k}=|A_{n,k}|=|A^H_{n,k}|+|A^U_{n,k}|=a_{n-1,k}+a_{n-1,k-1}.$$

\medskip

A related statistic that one can consider is the $x$-coordinate of the first descent in a bargraph. Let $\xfd$ denote this statistic on $\B$, and also the corresponding statistic on $\M$ under $\Delta$, which is the number of $H$ steps before the first descent. Then
$$M_{\xfd}=1+txM_{\xfd}+y(M_{\xfd}-1)+xy(M_{\xfd}-1)M,$$
and
$$B_{\xfd}=y(M_{\xfd}-1)=
ty\,\frac { 1+x-y-2tx-xy-\sqrt {(1-y) ( 1-2x-y-2xy+{x}^{2}-{x}^{2}y ) }  }{2(1-t-(1-t)y+({t
}^{2}-t)x+txy)}.$$

\subsection{Number of columns of height $h$ and number of initial columns of height $1$} \label{ch}

Let $\uc$ denote the number of columns of height 1 in bargraphs. The corresponding statistic on $\M$, via $\Delta$, is the number of $H$ steps on the $x$-axis, which we also denote by $\uc$.
Using the usual decomposition,
$$M_{\uc}=(1+y(M-1))(1+txM_{\uc})$$
and 
$$B_{\uc}=y(M_{\uc}-1)=y\,\frac {
1-x-y-xy+2t(1-t){x}^{2}(1-y)- \sqrt { (1-y) (1-2x-y-2xy+{x}^{2}-{x}^{2}y ) } }{2x ( 1-t+ty+({t}^{2}-t)x(1-y) )}.$$

We can now recursively count the number of columns of height $h$ in bargraphs, which we denote by $\ch_h$, for any $h\ge2$. They correspond to $H$ steps at height $h-1$ in cornerless Motzkin paths. From the decomposition,
$$M_{\ch_h}=(1+y(M_{\ch_{h-1}}-1))(1+xM_{\ch_h}),$$
and so 
$$M_{\ch_h}=\frac{1}{\dfrac{1}{1+y(M_{\ch_{h-1}}-1)}-x},$$
from where one can find $B_{\ch_h}=y(M_{\ch_h}-1)$. The above recurrence, in terms of bargraphs, can be written as
$$B_{\ch_h}= \frac{y}{x} \left(\frac{1}{1-x-xB_{\ch_{h-1}}} - 1 - x\right).$$

A related statistic is the number of initial columns of height 1, which we denote by $\iuc$. The corresponding statistic on $\M$ is the number of initial $H$ steps. We now have
$$M_{\iuc}=1+txM_{\iuc}+y(M-1)+xy(M-1)M$$
and
\begin{equation}\label{eq:Biuc}
B_{\iuc}=y(M_{\iuc}-1)=
\frac {1-2x-y+{x}^{2}+(2t-1){x}^{2}y
-(1-x)\sqrt {(1-y) [(1-x)^2-y(1+x)^2]}}{2x (1-tx ) }.
\end{equation}

\subsection{Least column height} \label{lch}

Let $\lch$ be the statistic ``least column height'' (that is, height of a shortest column) in a bargraph. The corresponding statistic on $\M$, which we also denote by $\lch$, is the least height of a horizontal step.
From the decomposition in the first line of Figure~\ref{fig:decomposition},
$$M_{\lch}=1+xM+ty(M_{\lch}-1)+xy(M-1)M,$$
since in the second and fourth cases there is a step at height 0, and in the third case the least height increases by one.
Solving for $M_{\lch}$ and using that $B_{\lch}=ty(M_{\lch}-1)$, we get
\begin{equation}\label{eq:Blch}
B_{\lch}=t(1-y)\,\frac {1-x-y-xy
-\sqrt {(1-y)[(1-x)^2-y(1+x)^2]}}{2x ( 1-ty)}.
\end{equation}

\subsection{Width of leftmost horizontal segment} \label{lhs}

Let $\lhs$ be the statistic ``width of leftmost horizontal segment,'' both in bargraphs and cornerless Motzkin paths. From the usual decomposition,
$$M_{\lhs}=1+\frac{tx(1-x)}{1-tx}M+y(M_{\lhs}-1)+xy(M_{\lhs}-1)M.$$
For the second summand, note that paths of the form $HA$ where $A\in\M$ starts with $H^\ell$ (but not $H^{\ell+1}$) contributes $t^{\ell+1}x^{\ell+1}(M-xM)$; summing over $\ell\ge0$ we get $\frac{tx(1-x)}{1-tx}M$.
Solving for $M_{\lhs}$ and using that $B_{\lhs}=y(M_{\lhs}-1)$, we get
\begin{equation}\label{eq:Blhs}
B_{\lhs}=t(1-x)\,\frac {1-x-y-xy
-\sqrt {(1-y)[(1-x)^2-y(1+x)^2]}}{2x ( 1-tx)}.
\end{equation}

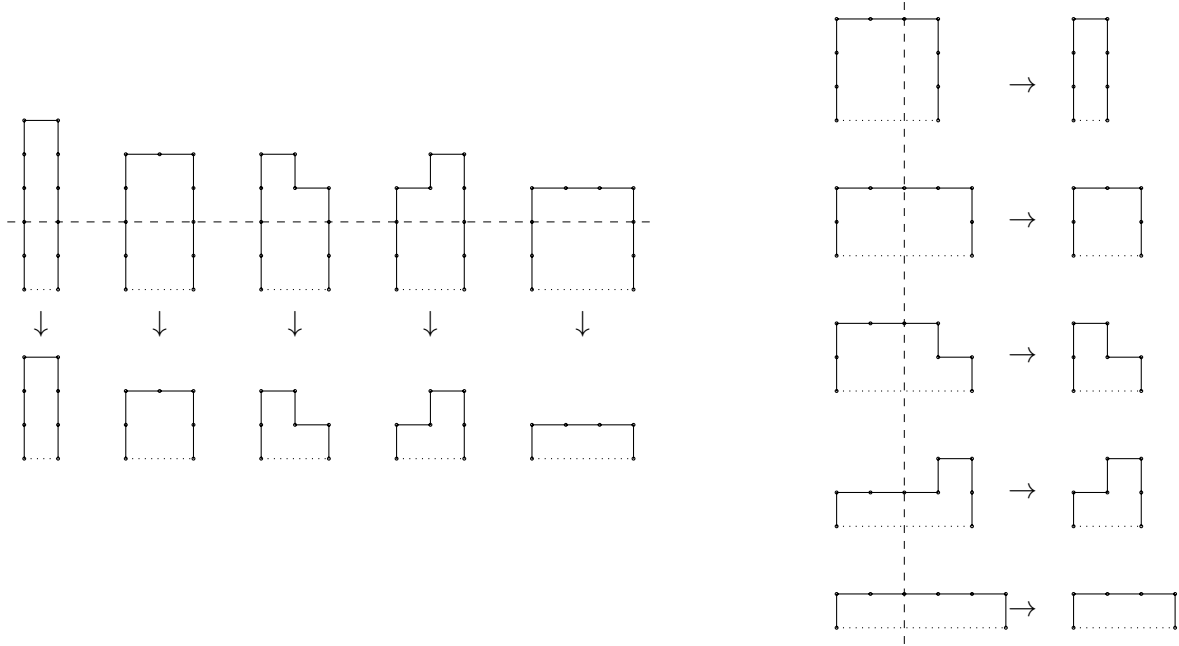
\begin{figure}[htb]
\centering
    \begin{tikzpicture}[scale=0.45]
     \draw(0,0) circle(1.2pt) \N\N\N\N\N\H\S\S\S\S\S;
     \draw(3,0) circle(1.2pt) \N\N\N\N\H\H\S\S\S\S;
     \draw(7,0) circle(1.2pt) \N\N\N\N\H\S\H\S\S\S;
     \draw(11,0) circle(1.2pt) \N\N\N\H\N\H\S\S\S\S;
     \draw(15,0) circle(1.2pt) \N\N\N\H\H\H\S\S\S;     
     
\draw[dotted] (0,0)--(1,0);
\draw[dotted] (3,0)--(5,0);
\draw[dotted] (7,0)--(9,0);
\draw[dotted] (11,0)--(13,0);
\draw[dotted] (15,0)--(18,0);
\draw[very thin,dashed] (-0.5,2)--(18.5,2);     
\draw (0.5,-1) node{$\downarrow$};
\draw (4,-1) node{$\downarrow$};
\draw (8,-1) node{$\downarrow$};
\draw (12,-1) node{$\downarrow$};
\draw (16.5,-1) node{$\downarrow$};

 \draw(0,-5) circle(1.2pt) \N\N\N\H\S\S\S;
     \draw(3,-5) circle(1.2pt) \N\N\H\H\S\S;
     \draw(7,-5) circle(1.2pt) \N\N\H\S\H\S;
     \draw(11,-5) circle(1.2pt) \N\H\N\H\S\S;
     \draw(15,-5) circle(1.2pt) \N\H\H\H\S;  
     
\draw[dotted] (0,-5)--(1,-5);
\draw[dotted] (3,-5)--(5,-5);
\draw[dotted] (7,-5)--(9,-5);
\draw[dotted] (11,-5)--(13,-5);
\draw[dotted] (15,-5)--(18,-5);     
     
\begin{scope}[shift={(24,-10)}]
     
\draw(0,15) circle(1.2pt) \N\N\N\H\H\H\S\S\S;
     \draw(0,11) circle(1.2pt) \N\N\H\H\H\H\S\S;
     \draw(0,7) circle(1.2pt) \N\N\H\H\H\S\H\S;
     \draw(0,3) circle(1.2pt) \N\H\H\H\N\H\S\S;
     \draw(0,0) circle(1.2pt) \N\H\H\H\H\H\S;   
\draw[dotted] (0,15)--(3,15);
\draw[dotted] (0,11)--(4,11);
\draw[dotted] (0,7)--(4,7);
\draw[dotted] (0,3)--(4,3);
\draw[dotted] (0,0)--(5,0);      

     \draw[very thin,dashed] (2,18.5)--(2,-0.5); 
     
\draw (5.5,16) node{$\rightarrow$};     
\draw (5.5,12) node{$\rightarrow$};     
\draw (5.5,8) node{$\rightarrow$};     
\draw (5.5,4) node{$\rightarrow$};     
\draw (5.5,0.5) node{$\rightarrow$};     
     
\draw(7,15) circle(1.2pt) \N\N\N\H\S\S\S;
     \draw(7,11) circle(1.2pt) \N\N\H\H\S\S;
     \draw(7,7) circle(1.2pt) \N\N\H\S\H\S;
     \draw(7,3) circle(1.2pt) \N\H\N\H\S\S;
     \draw(7,0) circle(1.2pt) \N\H\H\H\S;   
\draw[dotted] (7,15)--(8,15);
\draw[dotted] (7,11)--(9,11);
\draw[dotted] (7,7)--(9,7);
\draw[dotted] (7,3)--(9,3);
\draw[dotted] (7,0)--(10,0);    
     
              \end{scope}  
    \end{tikzpicture}
   \caption{The bijections $\{A\in\B_n:\lch(A)>h\}\to\B_{n-h}$ (left) and $\{A\in\B_n:\lhs(A)>h\}\to\B_{n-h}$ (right) for $n=6$ and $h=2$.}\label{fig:lch_lhs}
\end{figure}

By comparing~\eqref{eq:Blch} and~\eqref{eq:Blhs}, we see that 
$(1 - x)(1 - ty)B_{\lch} = (1 - y)(1 - tx)B_{\lhs}$. In particular,
$$B_{\lch}(t,z,z)=B_{\lhs}(t,z,z),$$
that is, the statistics $\lch$ and $\lhs$ are equidistributed on bargraphs of fixed semiperimeter.
Let us describe a direct combinatorial proof of this equality by showing that both
$\{A\in\B_n:\lch(A)>h\}$ and $\{A\in\B_n:\lhs(A)>h\}$ are in bijection with $\B_{n-h}$.
The bijection $\{A\in\B_n:\lch(A)>h\}\to\B_{n-h}$ is obtained by deleting the initial $U^h$ and the final $D^h$.
The bijection $\{A\in\B_n:\lhs(A)>h\}\to\B_{n-h}$ is obtained by deleting the initial $H^h$ from the leftmost horizontal segment. Examples of these bijections for all bargraphs with $n=6$ and $h=2$ are given in Figure~\ref{fig:lch_lhs}.

Another remarkable equidistribution phenomenon arises when comparing the statistics $\lhs$ and $\iuc$. In terms of generating functions, it can be checked by comparing~\eqref{eq:Biuc} and~\eqref{eq:Blhs} that
$$B_{\lhs}-\frac{txy}{1-tx}=t\left(B_{\iuc}-\frac{txy}{1-tx}\right).$$
Combinatorially, this formula states that if we ignore bargraphs consisting of one row ---which contribute $\frac{txy}{1-tx}$ to both generating functions---, then the distribution of the statistics $\lhs$ and $\iuc+1$ on the remaining bargraphs with a fixed number of $U$ and $H$ steps is the same. In particular, these statistics are equidistributed on bargraphs of fixed semiperimeter (and at least two rows), that is,
$$|\{A\in\B_n:\lhs(A)=h\}|=|\{A\in\B_n:\iuc(A)=h-1\}|$$ for all $1\le h\le n-2$.

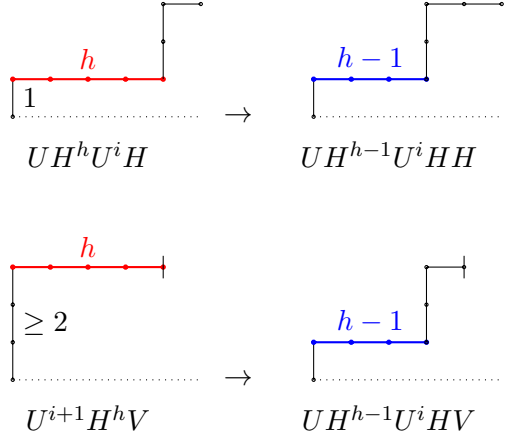
\begin{figure}[htb]
\centering
    \begin{tikzpicture}[scale=0.5]
     \draw(0,0) circle(1.2pt) \N; 
      \draw[thick,red](0,1)  circle(1.2pt) \H\H\H\H;
      \draw(4,1) \N\N\H; 
      \draw[red] (2,1) node[above] {$h$};
       \draw(0,0.5) node[right] {$1$};
      \draw[thin,dotted](0,0)--(5,0);
      \draw(2,-1) node {$UH^hU^iH$}; 
      \draw(6,0) node {$\rightarrow$};           
       \begin{scope}[shift={(8,0)}]
     \draw(0,0) circle(1.2pt) \N; 
      \draw[thick,blue](0,1)  circle(1.2pt) \H\H\H;
      \draw[blue] (1.5,1) node[above] {$h-1$};
       \draw(3,1)  circle(1.2pt) \N\N\H\H;
      \draw[thin,dotted](0,0)--(5,0);
      \draw(2,-1) node {$UH^{h-1}U^iHH$}; 
       \end{scope}
        \begin{scope}[shift={(0,-7)}]
      \draw(0,0) circle(1.2pt) \N\N\N; 
      \draw[thick,red](0,3)  circle(1.2pt) \H\H\H\H;
      \draw(4,3.3)--(4,2.7);
      \draw[red] (2,3) node[above] {$h$};
       \draw(0,1.5) node[right] {$\ge2$};
      \draw[thin,dotted](0,0)--(5,0);
      \draw(2,-1) node {$U^{i+1}H^hV$}; 
      \draw(6,0) node {$\rightarrow$};           
       \begin{scope}[shift={(8,0)}]
     \draw(0,0) circle(1.2pt) \N; 
      \draw[thick,blue](0,1)  circle(1.2pt) \H\H\H;
      \draw(4,3.3)--(4,2.7);
      \draw[blue] (1.5,1) node[above] {$h-1$};
       \draw(3,1)  circle(1.2pt) \N\N\H;
      \draw[thin,dotted](0,0)--(5,0);
      \draw(2,-1) node {$UH^{h-1}U^iHV$}; 
       \end{scope}
              \end{scope}  
    \end{tikzpicture}
   \caption{The bijection $\phi$.}\label{fig:phi}
\end{figure}

Next we give a combinatorial proof of this fact, by providing a bijection $\phi$ from 
the set of bargraphs with at least two rows to itself that preserves the number of $U$s and the number of $H$s, and such that $\iuc(\phi(A))=\lhs(A)-1$.
Given a bargraph A with at least two rows, let $h=\lhs(A)$. Consider two cases, as illustrated in Figure~\ref{fig:phi}:
\begin{itemize}
\item If the first column has height $1$, then $A$ starts with $UH^hU^iH$ (for some $i\ge1$). 
Replace this initial piece with $UH^{h-1}U^iHH$, and let $\phi(A)$ be the resulting bargraph. Viewing the bargraph as a polyomino, the height of column $h$ has been increased to match the height of column $h+1$.
\item If the first column has height at least $2$, then $A$ starts with $U^{i+1}H^hV$ (for some $i\ge1$), with $V\in\{U,D\}$. Replace this initial piece with $UH^{h-1}U^iHV$, and let $\phi(A)$ be the resulting bargraph. In terms of the polyomino, the first $h-1$ columns have been shortened to height~$1$.
\end{itemize}
In both cases, it is clear that $\iuc(\phi(A))=h-1$. Note that in the second case, if $h=1$, then $\phi(A)=A$.

To see that $\phi$ is a bijection, observe that every bargraph $A'$ with at least two rows starts either with $UH^{h-1}U^iHH$ or with $UH^{h-1}U^iHV$ (where $V\in\{U,D\}$), for some $i,h\ge1$, where $\iuc(A')=h-1$. Indeed, these cases are determined by whether the first horizontal segment that is not at height 1 has length 1 or at least 2, respectively. These cases correspond to the two bullets above, and so $\phi^{-1}(A')$ is obtained by replacing  $UH^{h-1}U^iHH$ with $UH^hU^iH$, or $UH^{h-1}U^iHV$ with $U^{i+1}H^hV$, respectively.

\subsection{Number of occurrences of $UHU$} \label{uhu}

Let $\uhu$ count the number of occurrences of $UHU$ in $\B$ and $\M$. From the decomposition, we have
$$M_{\uhu}=[1+y\,(M_{\uhu}-1+(t-1)x(M_{\uhu}-1-xM_{\uhu}))]\,(1+xM_{\uhu}).$$
This is because when the elevated path in the first piece starts with a $HU$, it creates an additional occurrence of $UHU$. Since the generating function with respect to $\uhu$ for paths that start with $H$ is $xM_{\uhu}$, the one for paths that start with $U$ is $M_{\uhu}-1-xM_{\uhu}$, and so the one for paths that start with $HU$ is $x(M_{\uhu}-1-xM_{\uhu})$.

Now, since $B_{\uhu}=y\,(M_{\uhu}-1+(t-1)x(M_{\uhu}-1-xM_{\uhu}))$, we get
that $B_{uhu}$ satisfies
$$xB_{uhu}^2 - (1 - x - y - txy)B_{uhu} + xy = 0.$$ 

\subsection{Length of the initial staircase $UHUH\dots$} \label{stair}

Let $\stair$ be the statistic that measures the length of the longest initial sequence of the form $UHUH\dots$ in a bargraph.
The corresponding statistic on $\M$, also denoted $\stair$, is the length of the longest initial sequence of the form $HUHU\dots$.

Using the first line of the decomposition in Figure~\ref{fig:decomposition},
\begin{align*}
\M&=\epsilon\cup H\M \cup [U\times(\M\setminus\epsilon)\times D] \times(\epsilon\cup H\M)\\
&=\epsilon\cup H \cup HH\M \cup [HU\times(\M\setminus\epsilon)\times D]\times (\epsilon\cup H\M)
\cup [U\times(\M\setminus\epsilon)\times D] \times(\epsilon\cup H\M),
\end{align*}
where to obtain the second line we have applied the decomposition again to the $\M$ in the term $H\M$. It now follows that 
$$M_{\stair}=1+tx+tx^2M+[t^2xy(M_{\stair}-1)](1+xM)+y(M-1)(1+xM).$$
Solving for $M_{\stair}$ and using that $B_{\stair}=ty(M_{\stair}-1)$, we get
$$B_{\stair}=t\,\frac { R -(1-x+tx^2-t^2xy)\sqrt{S}}{2x ( 1-{t}^{2}x+{t}^{2}{x}^
{2}-{t}^{2}xy+({t}^{4}-t^2){x}^{2}y) },$$
where $R=1-2x-y+(1+t)x^2-{t}^{2}xy+({t}^{2}+t-1){x}^{2}y+{t}^{2}x{y}^{2}-t{x}^{3}+(t-2t^3){x}^{3}y+{t}^{2}{x}^{2}{y}^{2}$ and $S=(1-y)(1-2x-y+x^2-2xy-x^2y).$

\subsection{Number of odd-height and even-height columns} \label{oheh}
Let $\oh$ and $\eh$ denote the statistics that count the number of columns of odd and even height in bargraphs, respectively. The corresponding statistics on $\M$, which we also denote by $\oh$ and $\eh$, are the number of $H$ steps at even and odd height, respectively. Let $B_{\oh,\eh}(t,s,x,y)$ and $M_{\oh,\eh}(t,s,x,y)$ be the generating functions on $\B$ and $\M$ where $t$ marks $\oh$ and $s$ marks $\eh$. Let $M_{\eh,\oh}(t,s,x,y)=M_{\oh,\eh}(s,t,x,y)$.

From the usual decomposition,
\begin{align*}
M_{\oh,\eh}&=(1+y(M_{\eh,\oh}-1))(1+txM_{\oh,\eh}),\\
M_{\eh,\oh}&=(1+y(M_{\oh,\eh}-1))(1+syM_{\eh,\oh}).
\end{align*}
The second equation can be obtained from the first one by interchanging the variables $t$ and $s$.
Solving for $M_{\oh,\eh}$ and using that $B_{\oh,\eh}=y(M_{\oh,\eh}-1)$,
we get
$$B_{\oh,\eh}=\frac {
(1-tx)(1-sy)-(1+tx)(1+sy){y}^{2}
-\sqrt {S}
}{2y ( s+t(1-s)x+tsxy) },
$$
where 
$$S=(1-y) ( 1-tx+y+txy ) 
 (1-tx-2sy+2t(s-1)xy+(s^2-1){y}^{2}-t({s}^{2}+1)x{y}^{2}-2s{y}^{3}+2ts(s-1)x{y}^{3}-{s}^{2}{y}^{4}-{s}^{2}tx{y}^{4}).$$

\subsection{Area}\label{area}
The enumeration of bargraphs according to semiperimeter and area was obtained in~\cite{BMB,BBK_parameters} by using a construction that adds one column at a time. Here we show that our decomposition of cornerless Motzkin paths can be used to derive a continued fraction. By $\area$ of a bargraph we mean the area of the region under the bargraph and above the $x$-axis. On cornerless Motzkin paths, the statistic $\area$ will refer to the area of the region directly below the $H$ steps of the path and above the $x$-axis (that is, the regions below $U$ and $D$ steps are not included). With this definition, we have 
\begin{equation}\label{eq:Barea} B_{\area}(t,x,y)=y(M_{\area}(t,tx,y)-1),
\end{equation}
since $\Delta$ lifts all the $H$ steps, and so each $H$ step contributes an additional unit to the area.

From the second line of the decomposition in Figure~\ref{fig:decomposition},
$$M_{\area}(t,x,y)=(1+y(M_{\area}(t,tx,y)-1))(1+xM_{\area}(t,x,y)).$$
Solving for $M_{\area}(t,x,y)$, we get
$$M_{\area}(t,x,y)=\frac{1}{\dfrac{1}{1+y(M_{\area}(t,tx,y)-1)}-x}.$$
Iterating this formula to get a continued fraction and using~\eqref{eq:Barea}, we get
$$B_{\area}=-y+\frac{y}{-tx+\dfrac{1}{1-y+\dfrac{y}{-t^2x+\dfrac{1}{1-y+\dfrac{y}{-t^3x+\dfrac{1}{\ddots}}}}}}.$$
Note that truncating the above continued fraction by stopping at $-t^nx+1$ produces the correct generating polynomial for bargraphs with up to $n$ vertical steps, since the contribution of the truncated terms contains a factor $y^{n+1}$.

\section{Subsets of bargraphs}\label{sec:subsets}

In this section we find generating functions for three subsets of bargraphs: symmetric, weakly alternating, and strictly alternating.

\subsection{Symmetric bargraphs}\label{sec:sym}

{\em Symmetric bargraphs} are those that are invariant under reflection along a vertical line. Let $\B^S$ denote the set of symmetric bargraphs.  These bargraphs correspond via the bijection $\Delta$ to paths in $\M$ that are symmetric, which we denote by $\M^S$. For convenience, the set $\M^S$ does not include the empty path. We use the term
{\it Motzkin path prefix} to refer to a path with steps $U$, $H$ and $D$ starting at the origin and not going under the $x$-axis, but ending at any height. Let $\P_h$ denote the set of Motzkin path prefixes without peaks and valleys ending at height $h$, and let $P_h=P_h(x,y)$ be corresponding generating function where $x$ marks the number of $H$s and $y$ marks the number of $U$s, as usual. Clearly, $P_0=M$. For $h\ge1$, we obtain
$$P_h=y P_{h-1}+x P_h + xy (M-1) P_h$$
by separating the cases when the path starts with a $U$ and does not return to the $x$-axis, 
when it starts with an $H$, and when starts with a $U$ and returns to the $x$-axis. In the latter case, the first return gives a decomposition of the form $UADHB$, where $A\in\M$ is non-empty and $B\in\P_h$. Solving for $P_h$, we get
$$P_h=\frac{yP_{h-1}}{1-x-xy(M-1)},$$
and iterating over $h$,
$$P_h=\frac{y^h M}{\left(1-x-xy(M-1)\right)^h}$$
for all $h\ge0$.

In a path in $\M^S$ with an odd number of steps, the middle step must be an $H$. 
In a path in $\M^S$ with an even number of steps, the two steps in the middle must be
$HH$, since otherwise it would contain a peak or a valley. It follows that every path in $\M^S$ is of the form $AHA'$ or $AHHA'$, where $A\in\P_h$ for some $h\ge0$, and $A'$ is the vertical reflection of $A$.
In that case, 
the number of horizontal steps in $AHA'$ or $AHHA'$ equals 
twice the number of horizontal steps in A plus $1$ or $2$, respectively, 
while the number of up steps in $AHA'$ or $AHHA'$ equals 
the number of up steps in $A$ plus the number of down steps in $A$, i.e. twice the 
number of up steps in $A$ minus $h$. It follows that
$$M^S(x,y)=\sum_{h\ge0}(x+x^2)\frac{P_h(x^2,y^2)}{y^h}=\sum_{h\ge0} 
\frac{(x+x^2) y^h M(x^2,y^2)}{\left(1-x^2-x^2y^2(M(x^2,y^2)-1)\right)^h}
=\frac{(x+x^2)M(x^2,y^2)}{1-\frac{y}{1-x^2-x^2y^2(M(x^2,y^2)-1)}}.
$$
The generating function for symmetric bargraphs is then $B^S=y M^S(x,y)$, and we get the expression
$$B^S=(1+x)\,
\frac { \sqrt { (1-y^2) [(1-x^2)^2- y^2(1+x^2)^2]}-1+{x}^{2}+{y}^{2}+2{x}^{2}y+{x}^{2}{y}^{2} }{ 2x\left(1-y-{x}^{2}-{x}^{2}y \right)}.$$

\subsection{Alternating bargraphs}\label{sec:alt}

{\em Weakly alternating bargraphs} are those where ascents and descents alternate. Equivalently, they are of the form
$$U^{i_1}H^{j_1}D^{k_1}H^{\ell_1}U^{i_2}H^{j_2}D^{k_2}H^{\ell_2}\dots U^{i_m}H^{j_m}D^{k_m}$$
for some $m\ge1$ and $i_r,j_r,k_r,\ell_r\ge1$ for all $r$.
{\em Strictly alternating bargraphs} (called {\em alternating bargraphs} in~\cite{PB,R_thesis}), are those of the form
$$U^{i_1}HD^{k_1}HU^{i_2}HD^{k_2}H\dots U^{i_m}HD^{k_m}$$
for some $m\ge1$ and $i_r,k_r\ge1$ for all $r$.
Here we derive generating functions for these sets of bargraphs by using our decomposition of cornerless Motzkin paths.

Let $\B^\WA$ be the set of weakly alternating bargraphs, and let $\M^\WA$ denote the set of alternating paths in $\M$, which we define as those whose ascents and descents alternate. Adapting the decomposition in the second line of Figure~\ref{fig:decomposition} to alternating paths, we get
\begin{equation}\label{eq:decompMA}
\M^\WA=(\epsilon\cup U\M^\Lambda D)\times(\epsilon\cup H\M^\WA),
\end{equation}
where $\M^\Lambda$ is the subset of $\M^\WA$ consisting of nonempty paths that do not start or end with an $H$, or consist only of $H$s.

Since $xM^\WA$ is the generating function for paths in $\M^\WA$ that start (respectively, end) with an $H$, and 
$x+x^2M^\WA$ is the generating function for paths in $\M^\WA$ that start and end with an $H$, we have, by inclusion-exclusion, that
$$M^\WA-xM^\WA-xM^\WA+(x+x^2M^\WA)=(1-x)^2M^\WA+x$$
is the generating function for paths that neither start nor end with an $H$. Removing the empty path and adding paths consisting only of $H$s, it follows that
$$M^\Lambda=(1-x)^2M^\WA+x-1+\frac{x}{1-x}.$$

By Equation~\eqref{eq:decompMA}, we have
$$M^\WA=(1+y M^\Lambda)(1+xM^\WA),$$
which we can solve to obtain $M^\WA$ and $M^\Lambda$. Finally, using that $B^\WA=yM^\Lambda$, we get
$$B^\WA=\frac {1-2x-y+2xy+{x}^{2}
-\sqrt {\left((1 - x)^2 - y\right)\left((1 - x)^2 - y(1 - 2x)^2\right)}
}
{2x(1-x)}.$$

To obtain the generating function $B^\SA$ for the set $\B^\SA$ of strictly alternating bargraphs, we note that weakly alternating bargraphs can be obtained uniquely from strictly alternating bargraphs by replacing each $H$ with an arbitrary non-empty sequence of $H$ steps. Therefore, 
$B^\WA(x,y)=B^\SA(\frac{x}{1-x},y)$, or equivalently, $B^\SA(x,y)=B^\WA(\frac{x}{1+x},y)$, giving
\begin{equation}\label{eq:BSA}
B^\SA=\frac {1-y+x^2y-\sqrt {1-2y+y^2-2x^2y-2x^2y^2+x^4y^2}
}{2x}.
\end{equation}

Interestingly, this generating function happens to be closely related with the generating function $K$ of Motzkin paths with no occurrences of $UD$, $UU$ and $DD$, where $x$ marks the number of $U$ and $D$ steps, and $y$ marks the number of $H$ steps. Indeed, letting $\mathcal{K}$ be the set of such Motzkin paths, it is clear that every nonempty path in $\mathcal{K}$ can be decomposed uniquely either as $HA$ or as $UA'DA$, where $A\in\mathcal{K}$ is arbitrary, and $A'$ is  any path in $\mathcal{K}$ that starts and ends with an $H$. Thus, $K$ satisfies the equation
$$K=1+yK+x^2(y+y^2K)K,$$ where the factor $y+y^2K$ is the contribution of the piece $A'$.
Solving this equation we obtain $$K=\frac {1-y-x^2y-\sqrt {1-2y+y^2-2x^2y-2x^2y^2+x^4y^2}
}{2x^2y^2}=\frac{B^\SA-xy}{xy^2},$$ where the last equality is obtained by comparing with Equation~\eqref{eq:BSA}. The generating function for these paths by semiperimeter, $K(z,z)$, gives sequence A023432 in~\cite{OEIS}.

Next we give a recursive bijection $f:\B^\SA\setminus\{UHD\}\to\mathcal{K}$ that maps the statistics $\#H$ (i.e., number of $H$ steps) and $\#U$ in bargraphs to the statistics $\#U+\#D+1$ and $\#H+2$ in Motzkin paths, respectively. In particular, the bijection takes bargraphs of semiperimeter $n$ (for $n\ge3$) to paths with $n-3$ steps.
The base case consists of bargraphs with only one $H$ step, for which we define $f(U^aHD^a)=H^{a-2}$, where $a\ge2$.

\begin{figure}[htb]
\centering
    \begin{tikzpicture}[scale=0.5]
     \fill[cyan] (0,2)--(0,6)--(1,6)--(1,4)--(2,4)--(2,5)--(3,5)--(3,2)--(0,2); 
     \draw (3.5,3) node[below] {$H$};
\fill[pink] (4,2)--(4,6)--(5,6)--(5,5)--(6,5)--(6,6)--(7,6)--(7,3)--(8,3)--(8,5)--(9,5)--(9,2)--(4,2); 
          \draw(0,0) circle(1.2pt)  \N\N\N\N\N\N\H\S\S\H\N\H\S\S\H\N\N\N\H\S\H\N\H\S\S\S\H\N\N\H\S\S\S\S\S; 
      \draw[red](0,0) circle(1.2pt) \N\N;    \draw[red](9,2) circle(1.2pt) \S\S;
      \draw[thin,dotted](0,0)--(9,0);
      \draw(11,1) node {$\rightarrow$};
      \draw(12,0) circle (1.2pt)  \H\H\U\H\U\H\H\D\H\D\U\H\H\U\H\D\H\D\H;  
       \draw[cyan,thick](15,1) circle (1.2pt)  \H\U\H\H\D;  
      \draw[pink,thick](22,0) circle (1.2pt)  \U\H\H\U\H\D\H\D\H;  
      \draw[red](12,0) circle (1.2pt)  \H\H;  
       \draw[thin,dotted](12,0)--(31,0);
    \end{tikzpicture}
   \caption{A strictly alternating bargraph and the corresponding Motzkin path in $\mathcal{K}$ obtained by applying~$f$.}\label{fig:f}
\end{figure}
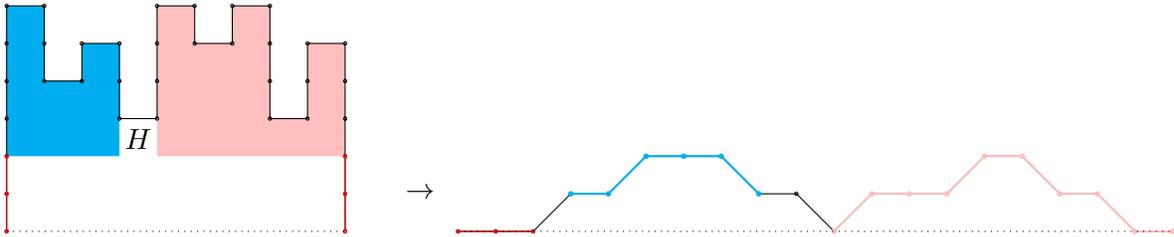

Every bargraph $G\in\B^\SA$ with at least two $H$ steps can be uniquely decomposed as 
\begin{equation}\label{eq:decompG} G=U^{a} G_1 H G_2 D^{a},
\end{equation} where $a\ge1$, $G_1\in\B^\SA$ and $UG_2D\in\B^\SA$.
Note that $a$ is the height of the lowest $H$ step in $G$, and that the $H$ that appears in the decomposition~\eqref{eq:decompG} is the leftmost $H$ step at height $a$. For such $G$, let
$$f(G)=H^{a-1} U\, f(UG_1D)\, HD\, f(UG_2D).$$
Note that $f(UG_1D)$ is either empty (when $G_1=UHD$) or it starts with an $H$. By induction, the path $f(G)$ contains no $UD$, $UU$ or $DD$. Figure~\ref{fig:f} shows an example of this construction.

\end{document}